\documentclass[11pt]{article}
\usepackage{mathrsfs}
\usepackage{latexsym,lineno}
\usepackage{epsfig}
\usepackage{color}
\usepackage{amsmath}\usepackage{verbatim}\usepackage{epsf}
\usepackage{amsthm}\usepackage{graphicx, float}\usepackage{graphicx}
\usepackage{amsfonts}\usepackage{amssymb}\usepackage{graphpap}
\usepackage{epic}\usepackage{curves}
\usepackage{tikz}
\usepackage{subfigure}
\newtheorem{fact}{Theorem}

\newtheorem{theorem}[fact]{Theorem}

\newtheorem{lemma}{Lemma}

\title{\bf Degree bipartite Ramsey numbers\thanks{Supported in part by NSFC and Shanghai Sailing Program(19YF1435500).} }
\author{Ye Wang$^{a}$, \hspace{2mm} Yusheng Li$^{b}$, \hspace{2mm} Yan Li$^{b}$\footnote{Corresponding author: 1610521@tongji.edu.cn} \\
{\small $^{a}$ School of Stats \& Math, \; Shanghai Lixin University of Accounting and Finance} \\
{\small Shanghai 201209, China}\\
{\small $^{b}$ School of Mathematical Sciences, Tongji University,\; Shanghai 200092, China} \\
{\small Email: wangye@sfu.edu.cn,\;li\_yusheng@tongji.edu.cn,\;1610521@tongji.edu.cn }   }
\date{}
\begin{document}
\date{}
\maketitle
\begin{abstract}
Let $H\xrightarrow{s} G$ denote that any edge-coloring of $H$ by $s$ colors contains a monochromatic $G$.
The degree Ramsey number $r_{\Delta}(G;s)$ is defined to be $\min\{\Delta(H):H\xrightarrow{s} G\}$,
and the degree bipartite Ramsey number $br_{\Delta}(G;s)$ is defined to be
$\min\{\Delta(H):H\xrightarrow{s} G\; \mbox{and} \;\chi(H)=2\}$.
In this note, we show that $r_{\Delta}(K_{m,n};s)$ is linear on $n$ with $m$ fixed.
We also determine $br_{\Delta}(G;s)$ where $G$ are trees, including stars and paths, and complete bipartite graphs.
\medskip

\end{abstract}

\section{Introduction}
\indent
For graphs $G$ and $H$, let $H\xrightarrow{s} G$ denote that any edge-coloring of $H$ by $s$ colors contains a monochromatic $G$.
The Ramsey number $r(G;s)$ is the smallest $N$ such that $K_N\xrightarrow{s} G$.
More generally, for any monotone graph parameter $\rho$, the $\rho$-Ramsey number is defined as
\[
r_\rho(G;s)=\min\{\rho(H):H\xrightarrow{s} G\}.
\]
This generalizes the Ramsey number since $r_\rho(G;s)=r(G;s)$
if $\rho(H)$ denotes the order of $H$.
When $\rho(H)$ denotes the size of $H$, it becomes the size Ramsey number $\hat{r}(G;s)$, see \cite{beck83,beck90,donadelli,erdos78,friedman,rodl}.
For cases $\rho(H)$ of being the clique number and the chromatic number of $H$,
we refer the reader, to the best of our knowledge, to \cite{folkman,lin,nesetril} and \cite{burr,zhu98,zhu11}, respectively.

The degree Ramsey number is defined as
\[
r_{\Delta}(G;s)=\min\{\Delta(H):H\xrightarrow{s} G\},
\]
where $\Delta(H)$ is the maximum degree of $H$.
Kinnersley, Milans and West \cite{kinnersley}, and Jiang, Milans and West \cite{jiang} obtained the bounds for degree Ramsey numbers of trees and cycles.
Kang and Perarnau \cite{kang} proved that $r_{\Delta}(C_4;s) = \Theta(s^2)$,
and Tait \cite{tait} proved that $r_{\Delta}(C_6;s) = \Theta(s^{3/2})$ and $r_{\Delta}(C_{10};s) = \Theta(s^{5/4})$.

In this note, let us define the degree bipartite Ramsey number $br_{\Delta}(G;s)$ as
\[
br_{\Delta}(G;s)=\min\{\Delta(H):\mbox{$H$ is bipartite and}\;\;H\xrightarrow{s} G\}.
\]
Obviously, for bipartite graph $G$, we have
\begin{equation}\label{e1}
r_{\Delta}(G;s)\leq br_{\Delta}(G;s).
\end{equation}
The results in this note reveal that the both sides of (\ref{e1}) are closed.
Now we consider the degree bipartite Ramsey numbers of trees, including stars, paths, and complete bipartite graphs.
Particularly, we show that $r_{\Delta}(K_{m,n};s)$ is linear on $n$ with $m$ fixed.
\begin{theorem}\label{t2}
If $T$ is a tree in which one vertex has degree $k$ and all others have degree at most $\lceil k/2\rceil$, then
$$br_{\Delta}(T;s)=s(k-1)+1.$$
\end{theorem}
Kinnersley, Milans and West \cite{kinnersley} showed $r_{\Delta}(T;s)=s(k-1)+1$
for any tree $T$ satisfying conditions in Theorem \ref{t2} with odd $k$, and thus the inequality in (\ref{e1}) is sharp.
They also proved
\begin{equation}\label{r-up}
r_{\Delta}(G;s)\le 2s(\Delta(T)-1).
\end{equation}
We shall generalize (\ref{r-up}) to the bipartite version.
Hence if the equality in (\ref{r-up}) holds, then the inequality in Theorem \ref{t3} becomes an equality from (\ref{e1}).
\begin{theorem}\label{t3}
If $T_n$ is a tree on $n$ vertices, then
$$br_{\Delta}(T_n;s)\leq 2s(\Delta(T)-1).$$
\end{theorem}

The above bound is sharp
since Alon, Ding, Oporowski and Vertigan \cite{alon} showed that
$r_{\Delta}(P_n;s)=2s$ for fixed $s$ and large $n$,
where $P_n$ is a path on $n$ vertices.

We now turn to the even cycles and complete bipartite graphs.

The results $r_\Delta(C_{2m};s)=\Theta(s^{1+\frac{1}{m-1}})$ in \cite{kang,tait}
and the upper bound in (\ref{e1}) imply $br_\Delta(C_{2m};s)\ge \Omega(s^{1+\frac{1}{m-1}})$
for cycles $C_{2m}$ with $m=2,3,5$ and $s\to \infty$.
It is well known that $ex(2n;C_{2m})\le O(n^{1+\frac{1}{m}})$
shown by Bondy and Simonovits \cite{bondy} for fixed $m$ and $n\to \infty$.
On the other hand, if the edges of $K_{n,n}$ are colored by $s$ colors,
then at least $n^2/s$ edges are monochromatic.
Therefore, if $n^2/s\ge ex(2n;C_{2m})$,
sufficiently $n^{1-\frac{1}{m}}\ge \Omega(s)$,
i.e., $\Delta(K_{n,n})=n\le O(s^{1+\frac{1}{m-1}})$,
any edge-coloring of $K_{n,n}$ by $s$ colors must contain a monochromatic $C_{2m}$ for such $n$.
This fact and the definition of $br_\Delta(C_{2m};s)$ imply
$br_\Delta(C_{2m};s)\le O(s^{1+\frac{1}{m-1}})$
for $m=2,3,5$.
Combining with the lower bound as mentioned, we have
\[
br_\Delta(C_{2m};s)=\Theta(s^{1+\frac{1}{m-1}}),
\]
for $m=2,3,5$ and $s\to \infty$.

The following result differs from the result in \cite{tait} which pointed out $r_\Delta(K_{m,n};s)=\Theta(s^m)$ for fixed $m$ and $n$ with $n>(m-1)!$ and $s\to \infty$.
\begin{theorem}\label{t1}
For fixed $m$ and $s$, if $n$ is large, then
\[
c n\le r_\Delta(K_{m,n};s)\le br_\Delta(K_{m,n};s)\le C n,
\]
where $c=c(m,s)$ and $C=C(m,s)$ are positive constants independent of $n$.
\end{theorem}
\medskip
\section{Proofs of main results}
\indent

\begin{lemma}\label{l1}
For any integers $n,s\geq 2$,
$br_{\Delta}(K_{1,n};s)=s(n-1)+1.$
\end{lemma}
{\bf Proof.}~~Since $K_{1,s(n-1)+1}\xrightarrow{s} K_{1,n}$,
we have $br_{\Delta}(K_{1,n};s)\le s(n-1)+1$.
For the lower bound, for any bipartite graph $H$ with maximum degree $s(n-1)$, let $H'$ be an $s(n-1)$-regular bipartite supergraph of $H$.
By Hall's Theorem, $H'$ decomposes into 1-factors.
Taking each of $s$ color classes to be the union of $n-1$ of these 1-factors yields an edge-coloring of $H'$ by $s$ colors with degree $n-1$ in each color at each vertex.
\hfill $\square$

{\bf Proof of Theorem \ref{t2}.}~~For the lower bound, for any tree $T$ in which one vertex has degree $k$ and all others have degree at most $\lceil k/2\rceil$,
it is obtained that $K_{1,k}\subseteq T$ and $br_{\Delta}(K_{1,k};s)\leq br_{\Delta}(T;s)$, so $br_{\Delta}(T;s)\geq s(k-1)+1$ by Lemma \ref{l1}.

For the upper bound, let $H$ be a regular bipartite graph having degree $s(k-1)+1$ and girth more than $|V(T)|$.
In any edge-coloring of $H$ by $s$ colors, some color class forms a spanning subgraph $H_1$ with average degree more than $k-1$.
If $H_1$ has a subgraph with a vertex $u$ of degree at most $r-1$ with $r=\lceil k/2 \rceil$,
as $k-1\geq 2(r-1)$, graph $H_1\setminus \{u\}$ has average degree more than $k-1$.
Thus there must be a subgraph $H_2$ in $H$ with minimum degree at least $r$ and average degree more than $k-1$.
Then $H_2$ also has a vertex of degree at least $k$, denoted by $v$.
In such a graph $H_2$, we can construct $T$ from $v$ by adding vertices.
When we want to add a vertex from a current leaf, it has $r-1$ neighbors in $H_2$ that (by the girth condition) are not already in the tree,
and then we get the desired tree $T$, finishing the proof. \hfill $\square$

The following lemma is a well known fact, and we shall use it to prove Theorem \ref{t3}. Here we sketch the proof. For graph $H$ with average degree $d>0$,
when we delete the vertices of degrees less than $d/2$ repeatedly
if any,
then the resulting graphs have non-decreasing average degrees and minimum degrees.
\begin{lemma}\label{l2}
For positive integers $\delta$ and $d$ with $d\geq 2(\delta-1)$,
if graph $H$ has average degree at least $d$,
then $H$ contains a subgraph with minimum degree at least $\delta$ and average degree at least $d$.
\end{lemma}

{\bf Proof of Theorem \ref{t3}.}~~Let $r=\Delta(T)$. And we can construct a $2s(r-1)$-regular bipartite graph $H$ with girth more than $|V(T)|$ which is known to be possible in various ways.
See for instance \cite{erdos} for constructing a $2s(r-1)$-regular graph $G$ with girth more than $|V(T)|$,
then the direct product $G\times K_2$ is a $2s(r-1)$-regular bipartite graph with girth $g(G\times K_2)\geq g(G)\geq |V(T)|$.

Consider an edge-coloring of bipartite graph $H$ by $s$ colors,
and some color class yields a monochromatic spanning bipartite subgraph $H_1$ with average degree at least $2(r-1)$.
By Lemma \ref{l2}, $H_1$ contains a subgraph $H_2$ with minimum degree at least $r$.
First we choose a vertex from $V(H_2)$ as the root of tree and then add new neighbors.
When we want to add from the current leaf, it has $r-1$ neighbors in $H_2$ that (by the girth condition) are not already in the tree.
Thus, we have the desired monochromatic tree $T$. \hfill $\square$

The following lemma appeared in \cite{kang} firstly,
and then it was restated by Tait \cite{tait} in a more general way.
Before stating it, we need some notations.
For $v\in V(G)$, denote by $N_G(v)$ the set of all neighbors of $v$ in $G$.
For graphs $G$ and $H$, a homomorphism $\phi$ from $G$ to $H$ is an edge preserving mapping from $V(G)$ to $V(H)$.
A homomorphism from $G$ to $H$ is locally injective
if $N_G(v)$ is mapped to $N_H(\phi(v))$ injectively
for every $v\in V(G)$.
A graph is $L_G$-free if it does not contain any graph in $L_G$ as a subgraph,
where $L_G$ is the set of all graphs $H$ such that there is a locally injective homomorphism from $G$ to $H$.

\begin{lemma}\cite{kang,tait}\label{l5}
Let $G$ be a graph with at least one cycle, $H$ a graph of maximum degree $\Delta$, and $f$ a monotone non-decreasing function. If the edges of $K_N$ can be partitioned into $f(N)O(N^{1-\xi})$ $L_G$-free graphs for fixed $\xi>0$, then $H$ can be partitioned into $f(200\Delta)O(\Delta^{1-\xi})$ graphs which are $G$-free.
\end{lemma}
Denote by $e(G)$ the size of graph $G$. The following lemma comes from the double counting by K\"{o}v\'{a}ri, S\'{o}s and Tur\'{a}n \cite{kovari}.
\begin{lemma}\label{l6}\cite{kovari}
Suppose that $G$ is a subgraph of $K_{M,N}$ with $e(G)\ge Np$ and
\[
N\binom{p}{m}>(n-1)\binom{M}{m},
\]
then $G$ contains $K_{m,n}$.
\end{lemma}

{\bf Proof of Theorem \ref{t1}.}~~
We shall show that
\[
e^{-2}s^{\frac{mn-1}{m+n}}n\le r_\Delta(K_{m,n};s)\le br_\Delta(K_{m,n};s)\le e^{s^2-1}s^{m}n,
\]
for large $n$.

The assertion is obvious for $m=1$ by Lemma \ref{l1},
so we assume $m\ge 2$.
Showing the lower bound is equivalent to showing that
any graph of maximum degree $\Delta$ can be partitioned into $(e^2\Delta/n)^{(m+n)/(mn-1)}$ graphs which are $K_{m,n}$-free.
By Lemma \ref{l5}, it suffices to show that $K_N$ can be partitioned into $(e^2N/n)^{(m+n)/(mn-1)}$ graphs which are $K_{m,n}$-free. Let us consider a random edge-coloring of $K_N$ by $s$ colors such that each edge is colored independently with probability $1/s$. Let $p$ be the probability that there is a monochromatic $K_{m,n}$. Then $p\le s\binom{N}{m+n}\binom{m+n}{m}/s^{mn}$.
For $s=(e^2N/n)^{(m+n)/(mn-1)}$, we have
\begin{eqnarray*}
p&\le& \Big(\frac{eN}{m+n}\Big)^{m+n}\Big(\frac{e(m+n)}{m}\Big)^{m}(\frac{1}{s})^{mn-1}
= \frac{e^{2m+n}N^{m+n}}{{m^m(m+n)}^{n}s^{mn-1}} \\
&=& \frac{n^{m+n}}{m^m e^n(m+n)^n}
=\frac{n^m}{m^me^n}(1-\frac{m}{m+n})^n \le \frac{n^m}{m^me^n}e^{-\frac{mn}{m+n}}\\
&\le& \frac{n^m}{(\sqrt{e}m)^m e^n},
\end{eqnarray*}
which implies if $n$ is sufficiently large, the probability that there is a monochromatic $K_{m,n}$ is less than one. Hence, $K_N\not\xrightarrow{s} K_{m,n}$, and the desired lower bound follows.

Let us consider a complete bipartite graph $K_{M,N}$ on bipartition $(A,B)$
and an edge partition $(E_1,E_2,\ldots,E_s)$.
Without loss of generality, assume that $|E_1|\ge MN/s$.
Then by setting $p=M/s$ in Lemma \ref{l6},
the subgraph induced by $E_1$ contains $K_{m,n}$ if
\[
N\binom{M/s}{m}>(n-1)\binom{M}{m}.
\]
Thus if we set $N=\lfloor \binom{M}{m}n/\binom{M/s}{m}\rfloor$,
for all $M\ge sm$,
\[
br_\Delta(K_{m,n};s)\le N\le Cn,
\]
where
\[
C=\binom{M}{m}/\binom{M/s}{m}=s^m \frac{(M-1)(M-2)\ldots[M-(m-1)]}{(M-s)(M-2s)\ldots[M-s(m-1)]}.
\]
For $m\le s+1$, taking $M=(s+1)m$, we have
\begin{eqnarray*}
C&\le& s^m\Big(\frac{M-(m-1)}{M-s(m-1)}\Big)^{m-1}= s^m \Big(1+\frac{(s-1)(m-1)}{m+s}\Big)^{m-1}\\
&\le& s^m e^{\frac{(s-1)(m-1)^2}{m+s}} \le  s^m e^{s^2-1},
\end{eqnarray*}
where we use the facts that $1+x\le e^x$ and the maximum value of $\frac{(s-1)(m-1)^2}{m+s}$ attains at $m=s+1$.

For $m\ge s+2$, taking $M=(m-1)^2$, we have
\begin{eqnarray*}
C&\le& s^m\Big(\frac{M-(m-1)}{M-s(m-1)}\Big)^{m-1}= s^m \Big(1+\frac{s-1}{m-s-1}\Big)^{m-1}\\
&\le& s^m e^{\frac{(s-1)(m-1)}{m-s-1}}\le  s^m e^{s^2-1},
\end{eqnarray*}
where the maximum value of $\frac{(s-1)(m-1)}{m-s-1}$ attains at $m=s+2$.
Hence, we have
$
br_\Delta(K_{m,n};s)\le e^{s^2-1}s^{m}n
$
as claimed. \hfill $\square$

\end{document}